\documentclass[12pt]{article}
\usepackage{amsmath,amssymb,amsthm,enumerate}
\title{On a sparse random graph with minimum degree {three}: {L}ikely Posa's sets are
large.}

\author{Alan Frieze\thanks{Frieze's research
supported in part by NSF}\,\,\,\,and Boris Pittel\thanks{
Pittel's research supported in part by NSF Grant  DMS-0805996}
\hbox{ }\\
\small Carnegie Mellon University  and
Ohio State University\\[-0.8ex]
\small\texttt{af1p@andrew.cmu.edu and bgp@math.ohio-state.edu}}
\usepackage[active]{srcltx}
\usepackage[usenames]{color}\usepackage{graphicx}
\definecolor{brown}{cmyk}{0, 0.72, 1, 0.45}
\definecolor{grey}{gray}{0.5}

\newcommand{\ignore}[1]{}

\def\th{\theta}
\begin{document}
\maketitle
{\center\small Mathematics Subject Classifications: 05C80, 05C30, 34E05, 60C05}
{\center\small Keywords: random, sparse graphs, degrees, longest path, Posa sets}
\def\si{\par\smallskip\noindent}
\def\bi{\par\bigskip\noindent}
\def\pr{\text{ P\/}}
\def\ex{\text{E\/}}
\def\de{\delta}
\def\eps{\varepsilon}
\def\la{\lambda}
\def\a{\alpha}
\def\be{\beta}
\def\ga{\gamma}
\def\part{\partial}
\def\G11{G_{1,1}(n,m)}
\def\tag#1 {\eqno(#1)}
\def\Cal{\mathcal}
\def\var{\text{Var\/}}
\newcommand{\card}[1]{\left|#1\right|}
\newtheorem{Theorem}{Theorem}[section]
\newtheorem{Lemma}{Lemma}[section]
\newtheorem{Proposition}{Proposition}[section]
\newtheorem{Corollary}{Corollary}[section]
\newtheorem{Remark}{Remark}[section]
\numberwithin{equation}{section}

\begin{abstract}
We consider the likely size of the endpoint sets produced by Posa rotations, when applied to
a longest path in a random graph with $cn,\,c\geq 2.7$ edges that
is conditioned to have minimum degree at
least three.
\end{abstract}

\section{Introduction}

In the pioneering paper~\cite{eqref9} Erd\H os and R\'enyi asked how
large $m$, the number of edges, should be for the uniformly random
graph on $n$ vertices ($G(n,m)$) with high probability (whp)
to have a Hamilton cycle. The problem was vigorously attacked
by the various authors, see references in Bollob\'as~\cite{eqref7}; in particular, Koml\'os
and Szemer\'edi~\cite{KS1} showed that $m=n^{1+\eps}$ suffices. A critical breakthrough
was achieved by P\'osa~\cite{Po}; he showed that $m=c n\ln n$, $c>30$,
is enough. Qualitatively this is the best possible, since $m=\Theta(n \ln n)$
edges are needed for $G(n,m)$ to be connected whp,~\cite{eqref9}. Progressively
stronger extensions of P\'osa's result for $G(n,m)$ were achieved by Korshunov\cite{Kor},
Koml\'os and Szemer\'edi~\cite{KS2}, {Ajtai, Koml\'os and Szemer\'edi \cite{AKS}}, Bollob\'as~\cite{Bo2}, Bollob\'as, Fenner and
Frieze~\cite{BFF}, and Bollob\'as and Frieze~\cite{BF}. The
proofs frequently used a deceptively simple but surprisingly potent
``P\'osa's Lemma'' from~\cite{Po}.

Here is the Lemma.  Given a graph $G$ and a vertex $x_0$,
let $P=x_0x_1\dots x_h$ be a longest path from
$x_0$. If  $(x_h,x_i)\in E(G)$ for some $i<h-1$, then $P^\prime =x_0\dots x_ix_hx_{h-1}\dots
x_{i+1}$ is also a path in $G$ from $x_0$, of the same edge length $h$. In words, $P^\prime$ is
obtained for $P$ by rotation via the edge $(x_i,x_h)$. Let $S$ consist of $x_h$ and  the set of
endpoints  of all paths obtainable from $P$ through any number of rotations. Let $T$ be
the set of all outside neighbors of $S$ in $G$, $T=N(S)$ in short.
$T\subset \{x_0,\dots,x_h\}$ as $P$ is a longest
path. P\'osa's Lemma states that there are no edges between $S$ and
the vertices of $P\setminus (S\cup T)$.

P\'osa's  Lemma implies that $|T|< 2|S|$. Now if a graph $G$
is not too sparse, one may expect that the not-too-large vertex sets $A$  are sufficiently
expanding, so that $|N(A)|\ge 2|A|$. If that is the case, then it follows from P\'osa's Lemma that
$|S|$ has to be ``large'', $|S|\ge s(G)$.  Bollob\'as's next crucial
observation was that, if a path $P$ cannot be
extended via a sequence of rotations at either of its ends and there is no cycle of  length
$h(P)+1$, then there are at least $\binom{s(G)}{2}$ ``non-edges'' with a property: adding any
such non-edge to $E(G)$ creates a cycle out of a properly rotated $P$.
Using these fundamental properties of
P\'osa's sets, and also de la Vega's
theorem on whp existence of long paths in $G(n,p)$, $(np>4\ln 2)$,~\cite{dlV}, Bollob\'as~\cite{Bo2},
\cite{eqref7} found a surprisingly direct proof of Korshunov's,
Koml\'os-Szemer\'edi's result on a sharp threshold value of $m$ and $p$ for whp-Hamiltonicity
of $G(n,m)$ and $G(n,p)$.

Specifically, he proved that, for $p$ in question, whp $s(G(n,p))=\Theta(n)$, that is
the likely number of those  beneficial non-edges is quadratic in $n$.
Using this, he identified a sequence $G(n,p_0)\subset G(n,p_1)
\subset G(n,p_k)\subset G(n,p)$, $k=k(n)$,
such that de la Vega's result applies to $G(n,p_0)$ with room to spare, and for each $j$, with the
conditional probability $1-O(n^{-2})$, the length of the longest
path in $G(n,p_{j+1})$ strictly exceeds that in $G(n,p_j)$,
if the latter is non-Hamiltonian. That $G(n,p_k)$ is
whp Hamiltonian was then immediate.

Later we used a broadly similar argument to show that a random
graph with minimum degree at least $2$, $G_2(n,m)$ in short, whp
has a perfect matching iff $G_2(n,m)$ has no
isolated odd cycles, Frieze and Pittel~\cite{FP2}. A counterpart of P\'osa's Lemma in our
case was a lemma inspired by Gallai--Edmonds Structure Theorem and Edmonds' Matching algorithm.
This lemma allowed us to prove that, in absence of a perfect matching,
with high conditional probability
there exist $\Theta(n^2)$ non-edges each of which would increase the maximum
matching number.  And the place of de la Vega's algorithm was taken up by Karp-Sipser Matching
Greedy~\cite{KaSi}, analyzed in detail in our earlier paper~\cite{FP1}.

As the title shows, the core of Bollob\'as' paper~\cite{Bo2} was a proof that the uniformly random
$d$-regular graph on $[n]$, $G_d(n)$,  is whp Hamiltonian, if $d>10^7$.
Fenner and Frieze~\cite{FF} {independently
proved that $d>796$ suffices and then Frieze~\cite{F} came up with an
algorithmic proof of a better bound $d>85$}. It was
commonly believed that $d\ge 3$ suffices, and indeed Robinson and Wormald~\cite{eqreRW} settled this
conjecture affirmatively. Their non-algorithmic proof was based on a refined version of the
second order moment prompted by their discovery that, for $d\ge 3$, $\ex[X_n^2]=O(\ex^2[X_n])$,
$X_n$ being the number of Hamilton cycles. (The random graphs $G(n,p), G(n,m)$, with
$p,m$ in question, lack this remarkable property.)
Very recently, Bohman and Frieze~\cite{BF} proved
that another well-known random graph, $G_{d-\text{out}}(n)$ is whp Hamiltonian, if $d\ge 3$, in which
case the average vertex degree is asymptotic to $6$.

Now, {consider the} whp Hamiltonicity of $G_3(n,m)$, the random graph on $[n]$ with $m$ edges and
minimum degree $3$, at least.  Of course,  $m\ge 3n/2$, and we had better assume $m > 3n/2$, since
equality implies that $G_3(n,m)=G_3(n)$. From a more general result in
Bollob\'as, {Cooper,} Fenner and Frieze~\cite{BCFF} it follows that
whp $G_3(n,m)$ is Hamiltonian if $m>128n$.
Our ultimate goal is to push this bound down,  close to the best possible $m>1.5n$, and to
construct an algorithm that finds a Hamilton cycle in {$O_p\bigl(n^{1+o(1)}\bigr)$} running
time for $m/n$ in the arising range.

\begin{Theorem}\label{Thm} Introduce $f_k(x)=\sum_{j\ge k}x^j/j!$, a tail of the series for
$e^x$,  $x^*=4.789771\dots$,  a unique positive root of
$$
\frac{x^3 f_1(x)}{f_2(x)^2}=1,
$$
and
$$
a^*=\frac{x^*f_2(x^*)}{2f_3(x^*)}\approx 2.6616.
$$
If $m\ge (a^*+\eps)n$, $\eps>0$, i. e. the average vertex degree
exceeds $5.32$, then whp for each pair of P\'osa's sets $(S,T)$,
$$
|S| +|T|\ge n^{1-\delta_n},\quad \delta_n =(\ln\ln n)^{-1/2}.
$$
In words, the likely P\'osa's sets are, at least, almost linear in size.
\end{Theorem}
\begin{Remark} If the vertex degree range is an arbitrary $D\subseteq [3,\infty)$, $3\in D$,
then the above assertion continues to hold if we replace $f_3(x)$ with
$$
f_D(x )=\sum_{j\in D}\frac{x^j}{j!},
$$
and $f_2(x)$, $f_1(x)$ with $f_D^\prime(x)$ and $f_D^{\prime\prime}(x)$ respectively. For instance,
if $D=\{3,4\}$, then the likely P\'osa's sets are almost linear in size if $m\ge 17/9\approx 1.9$.
\end{Remark}

The proof of this claim takes up the rest of the paper. It is quite technical, apparently due to
exceeding sparseness of $G_3(n,m)$ for $m/n$ that close to the best $1.5$.
We firmly believe that, in a complete analogy with Bollob\'as' proof of Hamiltonicity of
$G(n,p)$ and $G(n,m)$, and our result on the existence of a perfect matching in $G_3(n,m)$,
the random graph $G_3(n,m)$ is whp Hamiltonian if $m\ge (a^*+\eps)n$. In fact,
we have already found a necessary modification of Karp-Sipser Matching Greedy. This
algorithm builds a 2-matching, viz. a spanning subgraph of $G_3(n,m)$ 
with maximum degree at most two.
A companion paper \cite{F2G} proves that whp this 
algorithm  produces a 2-matching $M$ with $O(\log n)$ components provided $c\geq 11$.
It then shows how to transform $M$ efficiently into a Hamilton cycle.

\section{Posa's sets in a graph with minimum degree $3$ at least}

Consider a graph $G$. Suppose that $P$ is a longest path
from a fixed vertex $v_0$. Posa's set is a pair $(S,T)$.
$S$ is the set
of endpoints of paths obtainable from $P$ through rotations via edges connecting an
endpoint and another vertex of a current path; $T=N(S)$,  the set of neighbors of $S$,
outside of $S$. The vertices of $T$ are all in $P$, since otherwise $P$ would not be a longest
path from $v_0$. Posa established a key property of $(S,T)$, namely
\begin{equation}\label{T2S}
|T| < 2|S|.
\end{equation}
Given a graph $G=(V,E)$ and $A\subseteq V$, we
use $G(A)$ to denote the subgraph of $G$ induced by $A$, and $e(A)$ to denote the
number of edges in $G(A)$.
\begin{Lemma}\label{STcomplex}
Suppose that the minimum degree of $G$ is $3$ at least.  Then
\begin{equation}\label{eST>}
e(S\cup T) > |S\cup T|,
\end{equation}
i. e. the edge density of $G(S\cup T)$ strictly exceeds $1$.
\end{Lemma}

{\bf Proof of Lemma \ref{STcomplex}\/} Let $Q$ be any path obtained from $P$ by
rotations. Posa observed that
\begin{equation}\label{posaobs}
\text{every $t\in T$ has an $S$-neighbor
on $Q$.}
\end{equation}
Introduce $T_1$, a set of all vertices $t\in T$ such that
$t$ has only one neighbor $s\in S$. Posa's observation
implies that, for every such pair $(s,t)$, $t$ and $s$ are neighbors on every path $Q$. In
particular, when $s$ is an endpoint, $t$ is next to $s$.
It follows then
that for any other  vertex $t^\prime\in T_1$ with a single neighbor
$s^\prime\in S$ we have $s^\prime\neq s$. Therefore $|T_1|\le |S|$.
\si

Let $D(S)$ denote the total degree of vertices in $S$, and let
$D_S(T)(\leq D(S))$ denote the total number of neighbors of $T$ in
$S$. Since each $t\in T\setminus T_1$ has at least two neighbors
in $S$, and $|T_1|\le |S|$, we have
\begin{equation}\label{DT>}
D_S(T)\ge |T_1| + 2(|T|-|T_1|)\ge 2|T|-|S|.
\end{equation}
Hence, as each $s\in S$ has degree $3$ at least,
\begin{equation}\label{DS+DTgeq}
D(S)+D_S(T)\ge 3|S|+2|T|-|S|=2(|S|+|T|).
\end{equation}
As
\begin{equation}\label{eST=}
2e(S\cup T)= D(S)+D_S(T) +2e(T),
\end{equation}
$e(T)$ being the number of edges in the subgraph $G(T)$  induced by $T$, we see that
$$
e(S\cup T)\ge |S\cup T|.
$$
The rest of the argument is needed to upgrade this to the strict inequality.
\si

From the proof of \eqref{DS+DTgeq},  and \eqref{eST=},
it follows that the edge density of $G(S\cup T)$ may be equal $1$ only if
\si
(a)
$$
D(S)+D_S(T)=2(|S|+|T|);
$$
(b) there are no edges in $T$;
\si
(c) $|T_1|=|S|$;
\si
 (d) each vertex in $S$ has exactly two neighbors in $S\cup (T\setminus T_1)$;
 \si
 (e) each vertex in $T_2:=T\setminus T_1$ has exactly two neighbors in $S$.
 \si
 If one of (c), (d), (e) is violated then $D(S)+D_S(T)>2(|S|+|T|)$.
\si

{\it Case\/} $T_2=\emptyset$. Given a path $P$, the $S$-vertices
are distributed over $P$ as subpaths of vertex length $i\ge 1$,
next neighbors of subpaths being  $T_1$-vertices.
Since $T_1$-vertices remain the neighbors of their single $S$-neighbors
on every path, and $|T_1|=|S|$,  there can be only paths of length $1$
and $2$, ``monomers'' and ``dimers''.
An endpoint $s$ of $P$ is a monomer, as its left neighbor is $t\in T_1$. There are no other
monomers in $P$, since an interior monomer would be flanked by two $T_1$-vertices,
sharing a common neighbor in $S$, which is impossible. Consequently,  the leftmost
subpath of $P$ is a dimer $s_1,s_2$, sandwiched between two vertices $t_1,t_2\in T_1$.
No rotation from $P$ can use either $t_i$, as $s_i$ is the only $S$-neighbor of $t_i$,
or $s_2$,  as $t_2\notin S$. If the rotation uses $s_1$, then
$t_1, s_1, s, t $ becomes the new leftmost dimer with $s_1$ retaining the left position. Of course,
if a rotation does not use $s_1$, then $t_1,s_1,s_2,t_2$
remains the leftmost subpath. So no sequence of rotations will make $s_1$ an endpoint. Contradiction.
\si

{\it Case\/} $T_2\neq\emptyset$. By (d)-(e),  the graph $G(S\cup T_2)$ is a disjoint union
of cycles. By (b), each cycle contains at least two
vertices from $S$. In fact, there is just one cycle,
since otherwise there would exist two vertices $s_1,s_2\in S$ such that no sequence of rotations
starting with a path with the endpoint $s_1$ would lead to a path ending at $s_2$.

Since there are no edges between vertices in $T_2$, two vertices from $T_2$ cannot be
neighbors on the cycle. And no two vertices from $S$ can be neighbors either. Otherwise,
there is an arc $s_1 s_2 t$, with $s_1,s_2\in S$, $t\in T_2$. Consider a path $Q$ that ends at
$s_1$.  We know that the left neighbor of $s_1$ in $Q$  is a $t_1\in T_1$.
By considering the rotation from $Q$ via the edge
$(s_1,s_2)$ we see that $s_2$ has another neighbor $s_3\in S$ distinct from $s_1$.
Hence $s_2$ has at least three neighbors in $S\cup T_2$, namely $s_1,s_3,t$. This
violates (d).
\si

Therefore the vertices from $S$ and from $T_2$ alternate on the cycle. Hence $|T_2|=|S|$,
whence
$$
|T|=|T_1|+|T_2| =2|S|,
$$
which violates Posa's inequality \eqref{T2S}.
\si
So the edge density of $G(S\cup T)$ exceeds $1$.\qed
\begin{Remark}\label{rem1}
The above argument needs to be refined if we want to put a bound on the time taken construct the end-point sets.
In this case suppose that we are doing a sequence of rotations with fixed endpoint $v_0$. If a rotation would produce
an endpoint that has been produced before in this sequence, then we do not do this rotation. This
limits the time spent producing endpoints, but it will reduce the number
of endpoints, but we will now argue that Lemma \ref{STcomplex} continues to hold. Indeed, all we have to observe is that
\eqref{posaobs} continues to hold. The argument being identical to Pos\'a's argument.
\end{Remark}

In the course of the proof, having assumed that the edge
density of $G(S\cup T)$ is $1$,  we saw  that then  $G(S\cup T)$
must be quite special. Namely $|T_1|=|S|$, and either
(1) $T_2=\emptyset$ and $G(S\cup T)$
is a cycle on $S$, with each of $T_1$-vertices attached to its own $S$-vertex, or
(2) $|T_2|=|S|$, and $G(S\cup T)$ is an alternating cycle on a bipartition $(S,T_2)$,
with  each of $S$-vertices hosting its own pendant vertex from  $T_1$.  The punch line was that
neither of these two graphs, each of edge density $1$,  can be a Posa's graph $G(S\cup T)$.
 \si

 In the next section we will show that in the random graph $G^{(3)}(n,m)$ whp
 no vertex subset $A$, with $|A|\le\eps_0 \ln n$  can induce a subgraph of edge density
 exceeding $1$. So, by Lemma \ref{STcomplex}, whp $|S|>\eps_0\ln n$. We will also
 show that whp the edge density of the induced subgraph is $1+o(1)$, for
every $A$, with $\eps_0 \ln n <|A|\le n^{1-o(1)}$. It is natural then to focus on the $o(n)$-Posa's
sets of edge density close to $1$, anticipating that the induced subgraphs $G(S\cup T)$
should interpolate between those two special, impossible, graphs. To prepare, let us have
a look at the deterministic properties of  $G(S\cup T)$ with an edge density close to $1$.
\si

Introduce $G^*=G^*(S\uplus (T\setminus T_1))$,  a subgraph on the vertex set
$S\cup (T\setminus T_1)$ whose edges have at least one end in $S$; so we
disregard edges of $G(S\cup T)$ between vertices of $T$, and also edges joining
the pendant vertices of $T_1$ to their respective $S$--``hosts''. For $v\in S\cup(T\setminus T_1)$,
let $d(v;G^*)$ denote the degree of $v$ in $G^*$; by the definition of $G^*$,
$\min_v d(v;G^*)\ge 2$. Introduce
$$
S_2=\{v\in S\,:\,d(v;G^*)=2\}.
$$
\newpage
 \begin{Lemma}\label{inter}\
 \begin{enumerate}[(i)]
 \item No vertex from $S_2$ can be a neighbor of both a
 vertex in $S_2$ and a vertex in $T\setminus T_1$.
 \item Suppose that
 \begin{equation}\label{defsigma}
 e(S\cup T) = (1+\sigma)(s+t),\quad s:=|S|,\, t:=|T|,
\end{equation}
for some $\sigma>0$. Then, denoting $|T_1|=t_1$,
\begin{align}
s-2\sigma(s+t)\le&\, t_1 \le\, s,\label{st1}\\
\sum_{v\in S\cup(T\setminus T_1)}\bigl[d(v;G^*)-2\bigr]\le&\, 2\sigma (s+t).\label{si+ti}
\end{align}
\end{enumerate}
\end{Lemma}
\noindent
{\it Remark.\/} Recalling that $t<2s$, the bound \eqref{st1} is not vacuous if  $\sigma\le 1/6$.
\si

{\bf Proof of Lemma \ref{inter}.\/} (i) Suppose that there are $s_1,s_2\in S_2$ and $t\in T\setminus
T_1$ such that $(s_1,s_2)$ and $(s_1,t)$ are edges in $G^*$. $s_2$ has a neighbor $t_1\in
T_1$,
since $s_2$'s degree in $G(S\cup T)$ is at least, whence exactly, $3$. Consider a path $P$ with
$s_2$ as its endpoint. $t_1$ is necessarily a penultimate vertex of $P$. Rotating $P$ via
the edge $(s_1,s_2)$ must make the right $P$-neighbor of $s_1$ a new endpoint. So $s_1$
has a neighbor in $S$ distinct from $s_2$, and $d(s_1;G^*)\ge 3$. Contradiction.
\si
(ii) First, using \eqref{DT>}, $D(S)\ge 3s$, and \eqref{eST=},
we obtain
$$
2(s+t)+(s-t_1)\le 2e(S\cup T)\le 2(1+\sigma)(s+t),
$$
which implies \eqref{st1} as $s-t_1\ge 0$. Second, the total vertex degree of $G^*$ is
$$
\sum_{v\in S\cup(T\setminus T_1)}d(v;G^*)=2e(S\cup T)-2t_1-2e(T).
$$
Therefore
\begin{align*}
\sum_{v\in S\cup(T\setminus T_1)}\bigl[d(v;G^*)-2\bigr]=&\,2e(S\cup T)-2t_1-2e(T)-
2(s+t-t_1)\\
=&\,2e(S\cup T)-2(s+t)-2e(T)\\
\le&\, 2(1+\sigma)(s+t) -2(s+t),
\end{align*}
which implies \eqref{si+ti}.\qed
\si

Introduce
\begin{equation}\label{deft2t3}
\begin{aligned}
S_3:=&\,S\setminus S_2,\\
T_2:=&\,\{v\in T\setminus T_1: d(v;G^*)=2\},\quad T_3:=(T\setminus T_1)\setminus T_2,
\end{aligned}
\end{equation}
and denote $s_i=|S_i|$, $t_i=|T_i|$; so $s=s_2+s_3$,
$t=t_1+t_2+t_3$. It follows from \eqref{si+ti} that
\begin{equation}\label{eS3T3}
\sum_{v\in S_3\cup T_3}\bigl[d(v;G^*)-2\bigr]\le 2\sigma(s+t),
\end{equation}
and then
\begin{equation}\label{s3+t3}
|S_3\cup T_3|=s_3+t_3\le 2\sigma(s+t).
\end{equation}

Let $\mu_1$ denote the total number of edges in the subgraph of $G^*(S\uplus(T\setminus T_1))$
induced by $S$, and let $\mu_2$ denote the total number of the remaining edges of
$G^*(S\uplus(T\setminus T_1))$, those joining vertices of $S$ and $T\setminus T_1$,
and set $\mu=\mu_1+\mu_2$. Clearly
\begin{align}
2\mu_1+\mu_2=&\,\sum_{v\in S}d(v;G^*),\label{2m1m2}\\
\mu_2=&\,\sum_{v\in T\setminus T_1}d(v;G^*).\label{m2}
\end{align}
Adding the equations \eqref{2m1m2} and \eqref{m2},
\begin{align}
\mu:=&\,\frac{1}{2}\sum_{v\in(S\cup T)\setminus T_1}d(v;G^*)
=\,s_2+t_2+\frac{1}{2}\sum_{v\in S_3\cup T_3}d(v;G^*)\notag\\
=&\,s+t-t_1+\frac{\xi_1+\xi_2}{2},\label{mu=}
\end{align}
where
\begin{equation}\label{xi's}
\xi_1:=\sum_{v\in S_3}\bigl[d(v;G^*)-2\bigr]\ge s_3,
\quad\xi_2:=\sum_{v\in T_3}\bigl[d(v;G^*)-2\bigr]\ge t_3.
\end{equation}
It follows from  \eqref{m2} and \eqref{2m1m2} that
\begin{equation}\label{mu1,mu2}
\begin{aligned}
\mu_1=&\,s-t+t_1+\frac{\xi_1-\xi_2}{2},\\
\mu_2=&\,2(t-t_1)+\xi_2.
\end{aligned}
\end{equation}
From \eqref{eS3T3},
\begin{equation}\label{xi<}
\xi_1+\xi_2\le 2\sigma(s+t).
\end{equation}
\qed

\begin{Remark}
We note here that Lemma \ref{inter} continues to hold under the restrictions described in Remark \ref{rem1}.
\end{Remark}

Let two disjoint sets, $S$ and $T$,
the partitions $S= S_2\cup S_3$,  $T=T_1\cup T_2\cup T_3$, and $\xi_1,\xi_2$ be given.
Let ${\cal N}(\bold S,\bold T,\boldsymbol\xi)$
denote the total number of the subgraphs
$G^*(S\uplus (T\setminus T_1))$, with $\mu_1,\mu_2$ determined by \eqref{mu1,mu2},
such that the constraints \eqref{st1}, whence the constraints \eqref{eS3T3}, \eqref{s3+t3}
and \eqref{xi<} hold for some $\sigma>0$.

\begin{Lemma}\label{calNSTxi}\
\begin{enumerate}[(i)]
 \item
\begin{equation}\label{NSTxi<}
\begin{aligned}
{\cal N}(\bold S,\bold T,\boldsymbol\xi) \le&\,
{\cal N}_1(\bold s,\bold t,\boldsymbol\xi),\\
{\cal N}_1(\bold s,\bold t,
\boldsymbol\xi):=&\,2^{-s_2-t_2-\mu_1}\,\frac{(2\mu_1+\mu_2)!}
{\mu_1!}\,\exp\left[O\bigl(\sigma (s+t)\bigr)\right].
\end{aligned}
\end{equation}
\item There exists  $\sigma_0\in (0,1)$ such that,  for $\sigma\le \sigma_0$ and
\begin{equation}\label{t-s}
(1+\sigma^{1/2})\,s\le t\le 2(1-\sigma^{1/2})\,s,
\end{equation}
a stronger bound holds:
\begin{multline}\label{NSTxibetter}
{\cal N}(\bold S,\bold T,\boldsymbol\xi) \le {\cal N}_2(\bold s,\bold t,\boldsymbol\xi):=
{\cal N}_1(\bold s,\bold t,\boldsymbol\xi)(s+t)^2\\
\times\exp\left[-(2s-t)\ln\frac{s}{2s-t}-(t-s)\ln\frac{s}{t-s}
+O\bigl(\sigma^{1/2} (s+t)\bigr)\right].
\end{multline}
\end{enumerate}
\end{Lemma}

{\bf Proof of Lemma \ref{calNSTxi}.\/}
It is well known, Bollob\'as~\cite{eqref5}, that   $g(\bold d)$,
the total number of graphs on $[\nu]$ with vertex degrees $d_1,\dots,d_{\nu}$,
and total vertex degree $2M:=\sum_i d_i$ satisfies
\begin{equation}\label{gboldd}
g(\bold d)\le (2M-1)!!\prod\limits_{i=1}^\nu\frac{1}{ d_i!}.
\end{equation}
Here is a bipartite counterpart of \eqref{gboldd}. Let $\nu_1$, $\nu_2$, and $\bold
{d}^\prime=(d_1^\prime,\dots, d_{\nu_1}^\prime)$,
$\bold{d}^{\prime\prime}=(d_1^{\prime\prime},\dots,d_{\nu_2}^{\prime\prime})$ be such that
$$
\sum_{i\in [\nu_1]}d_i^\prime=\sum_{j\in [\nu_2]}d_j^{\prime\prime}=M.
$$
Denote by $g(\bold{d}^\prime,\bold{d}^{\prime\prime})$ the total number of bipartite graphs on
a bipartition $[\nu_1]\uplus [\nu_2]$, with the left vertices
and the right vertices having degrees $\bold{d}^\prime$ and $\bold{d}^{\prime\prime}$.
Then
\begin{equation}\label{gd'd''}
g(\bold{d}^\prime,\bold{d}^{\prime\prime})\le M!\prod_{i\in [\nu_1]}
\frac{1}{d_i^\prime!}\prod_{j\in [\nu_2]}\frac{1}{d_j^{\prime\prime}!}.
\end{equation}

(i) Let $\bold d=\{d_v\}_{v\in S}$ be the (generic) vertex degrees of a subgraph of
$G^*(S\uplus (T\setminus T_1))$ induced by $S$;
 so
\begin{equation}\label{dsum}
\sum_{v\in S}d_v=2\mu_1.
\end{equation}
Let $\bold{d}^\prime=\{d_v^\prime\}_{v\in S}$ and $\bold{d}^{\prime\prime}
=\{d_v^{\prime\prime}\}_{v\in T_2\cup T_3}$ denote the vertex degrees of the complementary
bipartite graph on the bipartition $S\uplus (T_2\cup T_3)$;  so
\begin{equation}\label{sumd'}
\sum_{v\in S}d_v^\prime=\sum_{v\in T_2\cup T_3}d_v^{\prime\prime}=\mu_2.
\end{equation}
Here $\mu_1,\mu_2,\mu=\mu_1+\mu_2$ are given by \eqref{2m1m2}, \eqref{m2} and
\eqref{mu=}. In addition,
\begin{equation}\label{d+d'}
d_v+d_v^\prime\,\,\left\{\aligned
&=2,\quad&& v\in S_2,\\
&\ge 3,\quad&&v\in S_3,\endaligned\right.
\end{equation}
and
\begin{equation}\label{d''}
d_v\,\,\left\{\aligned
&=2,\quad&& v\in T_2,\\
&\ge 3,\quad&&v\in T_3.\endaligned\right.
\end{equation}
Using \eqref{gboldd} and \eqref{gd'd''}, we get an upper bound for the number of the graphs
with vertex degrees $\bold{d},\bold{d}^\prime,\bold{d}^{\prime\prime}$:
$$
(2\mu_1-1)!! \,\mu_2!\,2^{-s_2-t_2}\prod_{v\in S_3}\frac{1}{d_v!\,d_v^\prime!}
\,\prod_{v\in T_3}\frac{1}{d_v^{\prime\prime}!}.
$$
Introducing
$$
f_k(x)=\sum_{j\ge k}\frac{x^j}{j!},
$$
we have then
\begin{multline*}
\sum_{\bold{d},\bold{d}^\prime,\bold{d}^{\prime\prime}\text{ meet }\atop\eqref{dsum}-\eqref{d''}}
\prod_{v\in S}\frac{1}{d_v!\,d_v^\prime!}\,\prod_{v\in T_2\cup T_3}\frac{1}{d_v^{\prime\prime}!}\\
= 2^{-t_2}\,[x^{2\mu_1}y^{\mu_2}]\,\left(\sum_{d+d^\prime=2}\frac{x^d y^{d^\prime}}{d!\,d^\prime!}
\right)^{s_2}\left(\sum_{d+d^\prime\ge 3}\frac{x^d y^{d^\prime}}{d!\,d^\prime!}\right)^{s_3}\cdot
[z^{\mu_2-2t_2}]\left(\sum_{d^{\prime\prime}\ge 3}\frac{z^{d^{\prime\prime}}}{d^{\prime\prime}!}
\right)^{t_3}\\
=2^{-t_2}\,[x^{2\mu_1}y^{\mu_2}]
\,\left[\frac{(x+y)^2}{2}\right]^{s_2} \!\!\bigl[f_3(x+y)\bigr]^{s_3}\cdot [z^{\mu_2-2t_2}]\,
\bigl[f_3(z)\bigr]^{t_3}\\
=2^{-s_2-t_2}\binom{2\mu_1+\mu_2}{\mu_2}\,[\xi^{2\mu_1+\mu_2}] \, \xi^{2s_2}f_3(\xi)^{s_3}
\cdot [z^{\mu_2-2t_2}]\,\bigl[f_3(z)\bigr]^{t_3}\\
=2^{-s_2-t_2}\binom{2\mu_1+\mu_2}{\mu_2}\,[\xi^{2\mu_1+\mu_2-2s_2}]\,f_3(\xi)^{s_3}
\cdot  [z^{\mu_2-2t_2}]\,\bigl[f_3(z)\bigr]^{t_3}\\
\le  2^{-s_2-t_2}\binom{2\mu_1+\mu_2}{\mu_2}\,f_3(1)^{s_3+t_3}.
\end{multline*}
Thus
\begin{multline*}
{\cal N}(\bold S,\bold T,\boldsymbol\xi)\le (2\mu_1-1)!!
\,\mu_2!\,2^{-s_2-t_2}\binom{2\mu_1+\mu_2}{\mu_2}\,f_3(1)^{s_3+t_3}
\\
=2^{-s_2-t_2}\frac{(2\mu_1-1)!!}{(2\mu_1)!}\,(2\mu_1+\mu_2)!\,f_3(1)^{s_3+t_3}.
\end{multline*}
So, by \eqref{xi's} and \eqref{xi<},
\begin{equation}\label{mu1mu2bound}
{\cal N}(\bold S,\bold T,\boldsymbol\xi)\le {\cal N}_1(\bold s,\bold t,
\boldsymbol\xi)
:= 2^{-s_2-t_2-\mu_1} \frac{(2\mu_1+\mu_2)!}{\mu_1!}\,\exp\bigl[O(\sigma (s+t))\bigr].
\end{equation}
\si

(ii) Let $\{d_v\}_{v\in S_2}$, $\{d_v\}_{v\in S_3}$ be the vertex
degrees of the subgraphs of $G^*(S\uplus(T\setminus T_1))$ induced by $S_2$ and $S_3$
respectively. Let $\{\delta_v\}_{v\in S_2}$,
$\{\delta_v\}_{v\in S_3}$ denote the vertex degrees of a
bipartite graph induced by the bipartition $S_2\uplus S_3$. Finally, let
$\{d^\prime_v\}_{v\in S}$, $\{d^{\prime\prime}_v\}_{v\in T
\setminus T_1}$denote the vertex degrees of a
bipartite graph induced by the bipartition  $S\uplus (T\setminus T_1)$. By the definition,
\begin{align*}
&d_v+\delta_v+d^\prime_v=2,\,\,(v\in S_2),\quad d_v+\delta_v+d^\prime_v\ge 3,\,\,
(v\in S_3),\\
&d^{\prime\prime}_v=2,\,\,(v\in T_2),\quad d^{\prime\prime}_v\ge 3,\,\,(v\in T_3),
\end{align*}
and  by Lemma \ref{inter}
$$
d_v\cdot d^\prime_v= 0,\quad v\in S_2.\quad (!)
$$
Denote
\begin{align*}
&\sum_{v\in S_2}d_v=2\nu_2,\quad \sum_{v\in S_2}\delta_v=\nu_{2,3},\quad
\sum_{v\in S_2}d^\prime_v=\mu_{2,2},\\
&\sum_{v\in S_3}d_v=2\nu_3,\quad\sum_{v\in S_3}\delta_v=\nu_{3,2},\quad
\sum_{v\in S_3}d^\prime_v=\mu_{3,2};
\end{align*}
then $\nu_{2,3}=\nu_{3,2}$, and
\begin{equation}\label{mu1=}
\begin{aligned}
&2\nu_2+\nu_{2,3}+\mu_{2,2}=2|S_2|=2s_2,\\
&2\nu_3+\nu_{3,2}+\mu_{3,2}=2\mu_1+\mu_2-2s_2,\\
&\mu_{2,2}+\mu_{3,2}=\mu_2,\\
&\nu_2+\nu_{2,3}+\nu_3=\mu_1.
\end{aligned}
\end{equation}
Given the values of $\nu_2,\nu_3,\nu_{2,3}$, and $\mu_{2,2},\mu_{3,2}$,
the number of the corresponding subgraphs $G^*(S\uplus(T\setminus T_1))$
is bounded, as in part (i), by
\begin{multline}\label{3lines}
(2\nu_2-1)!!\,(2\nu_3-1)!!\,\nu_{2,3}!\,\mu_2! \,2^{-t_2}\,f_3(1)^{t_3}\\
\times \bigl[x_1^{{2}\nu_2}\,x_2^{\mu_{2,2}}\,x_3^{\nu_{2,3}}\bigr]\,
\left(\sum_{d+d^\prime+\delta=2\atop d\cdot d^\prime =0}\frac{x_1^d\,x_2^{d^\prime}\,
x_3^{\delta}}{d!\,d^\prime!\,\delta!}\right)^{s_2}\\
\times\bigl[y_1^{2\nu_3}\,y_2^{\mu_{3,2}}\,y_3^{\nu_{3,2}}\bigr]
\left(\sum_{d+d^\prime+\delta\ge 3}\frac{y_1^d\,y_2^{d^\prime}\,y_3^{\delta}}{d!\,d^\prime!\,\delta!}
\right)^{s_3}.
\end{multline}
The last line factor is
\begin{multline}
\binom{2\nu_3+\mu_{3,2}+\nu_{3,2}}{2\nu_3,\,\mu_{3,2},\,\nu_{3,2}}
\,\bigl[y^{2\nu_3+\mu_{3,2}+\nu_{3,2}}\bigr]\,f_3(y)^{s_3}\\
\le3^{2\nu_3+\mu_{3,2}+\nu_{3,2}}\,f_3(1)^{s_3}=
 3^{2\mu_1+\mu_2-2s_2}\,f_3(1)^{s_3}.\label{add1}
  \end{multline}
By \eqref{mu1,mu2} and \eqref{xi<}, this is
\begin{equation}\label{add2}
2\mu_1+\mu_2-2s_2=2s_3+\xi_1.
\end{equation}
Now $s_3\leq \xi_1$, see \eqref{xi's}, and so the RHS of \eqref{add1} is bounded by
\begin{equation}\label{third}
3^{2s_3+\xi_1}\, f_3(1)^{s_3}\le \bigl(27f_3(1)\bigr)^{\xi_1}\le \bigl(27f_3(1)\bigr)^{2\sigma(s+t)}.
\end{equation}
The second line factor in \eqref{3lines} is bounded by
\begin{equation}\label{secline}
\frac{\left(\frac{x_1^2}{2}+x_1x_3+
\frac{(x_2+x_3)^2}{2}\right)^{s_2}}{x_1^{{2}\nu_2}\,x_2^{\mu_{2,2}}\,x_3^{\nu_{2,3}}},
\end{equation}
for all $x_1,x_2,x_3>0$. The challenge is to select the ``best'' $x_1,x_2,x_3$. First of all
$\mu_{3,2},\,\nu_{2,3}\le 3\xi_1$, since $\nu_{2,3}=\nu_{3,2}$
and by \eqref{xi's}, \eqref{mu1=} and \eqref{add2}
$$
2\nu_3+\mu_{3,2}+\nu_{3,2}\le 3\xi_1.
$$
Therefore, by  \eqref{st1}, \eqref{eS3T3}, \eqref{mu1,mu2}, \eqref{xi<} and \eqref{mu1=},
\begin{align}
\nu_2=&\,\mu_1-\nu_{2,3}-\nu_3\leq 2s-t+O(\sigma(s+t)),\label{nu2=}\\
\mu_{2,2}=&\,\mu_2-\mu_{3,2}=2(t-s)+O(\sigma(t+s)).\label{mu22=}
\end{align}
By the condition in the Lemma,
the explicit terms are of order $\sigma^{1/2}(s+t)$ at least, thus dwarf
the remainders if $\sigma$ is small. We pick
$$
x_1=(2\nu_2)^{1/2},\quad x_2=(\mu_{2,2})^{1/2},\quad x_3= (\nu_{2,3})^{1/2}.
$$
Then
\begin{align*}
\frac{x_1^2}{2}+x_1x_3+\frac{(x_2+x_3)^2}{2}=&\,\frac{1}{2}(x_1^2+x_2^2+x_3^2)+(x_1+x_2)
x_3\\
=&s_2+O\bigl(\sqrt{(s+t)\xi}\bigr)\\
=&\,s+O\bigl(\sigma^{1/2}(s+t)\bigr),
\end{align*}
as $s_2=s+O\bigl(\sigma (s+t)\bigr)$. So the fraction in \eqref{secline} can be bounded from above by
\begin{align}
&\frac{1}{2^{s_2}}\sqrt{\frac{(2s_2)^{2s_2}}{(2\nu_2)^{2\nu_2}\,\mu_{2,2}^{\mu_{2,2}}\,\nu_{2,3}
^{\nu_{2,3}}}}\,\exp\bigl(O(\sigma^{1/2}(s+t))\bigr)\nonumber\\
&=\frac{1}{2^{s_2}}\exp\left(\nu_2\ln\frac{s_2}{\nu_2}+\frac{\mu_{2,2}}{2}\ln\frac{s_2}{\mu_{2,2}/2}
+\frac{\nu_{2,3}}{2}\ln\frac{s_2}{\nu_{2,3}/2}\right)\,\exp\bigl(O(\sigma^{1/2}(s+t))\bigr)\nonumber\\
&=\frac{1}{2^{s_2}}\exp\left[(2s-t)\ln\frac{s}{2s-t}+(t-s)\ln\frac{s}{t-s}+O\bigl(
\sigma^{1/2}(s+t)\bigr)\right].\label{second}
\end{align}
Turn to the first line of \eqref{3lines}. Using $(2a-1)!!\le 2^a a!$, and \eqref{mu1=}, \eqref{nu2=},
 we get
\begin{align}
(2\nu_2-1)!!\,(2\nu_3-1)!!\,\nu_{2,3}!\,\mu_2!
&\le\, 2^{\nu_2+\nu_3}\nu_2!\,\nu_{2,3}!\,\nu_3!\,\mu_2!
\le 2^{\nu_2+\nu_3}\mu_1!\,\mu_2!\nonumber\\
&=\,2^{\mu_1}\mu_1!\,\mu_2!\,\exp\bigl(O(\sigma(s+t))\bigr).\label{first}
\end{align}
Putting together \eqref{3lines}, \eqref{third}, \eqref{second} and \eqref{first}, we conclude
that the number of subgraphs $G^*(S\uplus(T\setminus T_1))$ with parameters $\boldsymbol
\mu$, $\boldsymbol\nu$ is bounded by
\begin{equation}\label{betterbound}
\frac{\mu_1!\,\mu_2!}{2^{s_2+t_2-\mu_1}}\exp\left[H(s,t)+O\bigl(\sigma^{1/2}(s+t)\bigr)\right],
\end{equation}
where
\begin{equation}\label{Hst}
H(s,t)=(2s-t)\ln\frac{s}{2s-t}+(t-s)\ln\frac{s}{t-s}.
\end{equation}
We emphasize that the remainder term estimate is uniform over the range of the parameters
$\boldsymbol\xi$, $\boldsymbol\nu$.

Let us compare the bounds \eqref{mu1mu2bound} and \eqref{betterbound}. We have
$$
\frac{2^{-s_2-t_2-\mu_1} \frac{(2\mu_1+\mu_2)!}{\mu_1!}}
{2^{-s_2-t_2+\mu_1}\,\mu_1!\,\mu_2!}=2^{-2\mu_1}\binom{2\mu_1+\mu_2}{\mu_1,\,\mu_1,\,\mu_2}.
$$
Using
$$
\mu_1=2s-t +O\bigl(\sigma(s+t)\bigr),\quad \mu_2=2(t-s)+O\bigl(\sigma(s+t)\bigr),
$$
(cf. \eqref{nu2=}, \eqref{mu22=}), and the Lemma condition on $2s-t$, $t-s$, it is simple to
show that the last expession is
$$
\exp\left[2H(s,t)+O\bigl(\sigma^{1/2}(s+t)\bigr)\right].
$$
Thus the bound \eqref{betterbound} can be written as
\begin{equation}\label{betterconnected}
{\cal N}_1(\bold s,\bold t)\exp\left[-2H(s,t)+O\bigl(\sigma^{1/2}(s+t)\bigr)\right].
\end{equation}
The bound \eqref{betterconnected} implies \eqref{NSTxibetter}, since
the factor $(s+t)^2$ is an upper
bound for the number of solutions $(\nu_2,\nu_{2,3},\nu_3,\mu_{2,2},\mu_{3,2})$ of \eqref{mu1=}.\qed
\si

Let $G^{**}(S\cup T)$ denote
$G^*(S\cup(T\setminus T_1))$ adorned with $t_1$ pendant $T_1$-vertices attached to
some $t_1$ $S$-vertices. That is, $G^{**}(S\cup T)$ is $G(S\cup T)$ without the edges joining
$T$-vertices to each other.

Given two disjoint sets $S$ and $T$,  let $N(\bold s,\bold t, \boldsymbol{\xi})$ denote
the total number of graphs $G^{**}(S\cup T)$ with $|S_2|=s_2$, $|T_1|=t_1$, $|T_2|=t_2$ and
the parameters $\xi_1,\xi_2$. The number of ways to select $T_1\subset T$ of cardinality
$t_1$ and then to match the vertices of $T_1$ with some $t_1$ vertices in $S$ is
 $t_1!\binom{s}{t_1}\binom{t}{t_1}$. The number of ways to select
$S_2\subset S$ of cardinality $s_2$ and to select $T_2\subset T\setminus T_1$ of cardinality
$t_2$ is $\binom{s}{s_2}\binom{t-t_1}{t_2}$, at most. (We neglect the constraint that $S_2$
needs to be a subset of the set of $S$-partners of $T_1$-vertices.)
The total count of possibilities is bounded by the product of those two.
\begin{Lemma}\label{Nu} In the notations of Lemma \ref{calNSTxi},
\begin{enumerate}[(i)]
\item
\begin{equation}\label{1Nu<}
N(\bold s,\bold t)\le t_1!\,\binom{s}{t_1}\binom{t}{t_1}\binom{s}{s_2}
\binom{t-t_1}{t_2}\,{\cal N}_1(\bold s,\bold t,\boldsymbol\xi);
\end{equation}
\item if $\sigma$ is small enough, and
$$
(1+\sigma^{1/2})\,s\le t\le 2(1-\sigma^{1/2})\,s,
$$
then
\begin{equation}\label{2Nu<}
N(\bold s,\bold t,\boldsymbol\xi)\le t_1!\,\binom{s}{t_1}\binom{t}{t_1}\binom{s}{s_2}
\binom{t-t_1}{t_2}\,{\cal N}_2(\bold s,\bold t,\boldsymbol\xi).
\end{equation}
\end{enumerate}
\end{Lemma}
\qed

Motivated by Lemma \ref{STcomplex} and Lemma \ref{inter},
in the next section we will focus on  subgraphs of $G^{(3)}(m,n)$ of size not exceeding $n^{1-o(1)}$,
showing as promised above that the likely edge density of such subgraphs is asymptotic to
$1$. It will remain to prove in the last section that whp  there are no Posa's
sets $S,T$ of low edge density, with $s+t=O\bigl(n^{1-o(1)}\bigr)$. Lemma \ref{Nu} will be
a key ingredient of that argument.

\section{Edge density of subgraphs of $G^{(3)}(n,m)$.}

Recalling the notation
$$
f_k(x)=\sum_{j\ge k}\frac{x^k}{k!},
$$
introduce $\lambda$,  a unique positive root of
\begin{equation}\label{la}
\frac{\la f_2(\la)}{f_3(\la)}=c:=\frac{2m}{n},
\end{equation}
\begin{Lemma}\label{e<(1+eps)v}\
\begin{enumerate}[(i)]
 \item
 For $\eps_0=(1/3)\ln ^{-1}\bigl[27f_3(8\la)/2^{10}\la f_2(\la)\bigr]$,
\begin{equation}\label{Pev}
\pr\bigl(\exists\,A\subset [n],\,|A|\leq \eps_0\ln n\,:\,e(A)>|A|\bigr)\to 0.
\end{equation}
Consequently, whp there does not exist an endpoint set $S$ of size below $\eps_0 \ln n$.
\item Let $\sigma_n\to 0$, but $(\sigma_n \ln n)/\ln\ln n\to\infty$,
and let $\rho_n=(\sigma_n \ln n)^{-1/2}$,
so that $\rho_n\to 0$. Then
\begin{equation}\label{density=1}
\pr\bigl(\exists\,A\subset [n],\,\eps_0\ln n\le |A|\leq n^{1-\rho_n}\,:\,
e(A)\ge (1+\sigma_n)|A|\bigr)\to 0.
\end{equation}
In words, with high probability, the edge density of subgraphs induced by sets $A$, of size
from $\eps_0 \ln n$ to $n^{1-o(1)}$,  is $1+o(1)$ at most.
\end{enumerate}
\end{Lemma}
\si

{\bf Proof of Lemma \ref{e<(1+eps)v}.\/} Our random
graph $G^{(3)}(n,m)$ is distributed uniformly on the set of
all $C(n,m)$ graphs of minimum degree at least $3$, with $m$ edges and $n$ vertices.
From a more general result in Pittel and Wormald~\cite{eqref21},  for $c>3/2$ we have:
for $n\to\infty$,
\begin{equation}\label{Cnm=}
C(n,m)\sim (2\pi n\var [Z])^{-1/2}(2m-1)!!\frac{f_3(\la)^n}{\la^{2m}}\,
\exp\bigl(-\eta -\eta^2/2\bigr);
\end{equation}
here $\la$ is the root of \eqref{la},  and
$Z$ is Poisson$(\la)$ conditioned on being at least 3. Probabilistically, \eqref{la} says that
$\ex[Z]=2m/n$. Also $\eta:=c^{-1}\ex\bigl[\binom{Z}{2}\bigr]$. Constant factors aside, the
claim is that
\begin{equation}\label{CTheta}
C(n,m)=\Theta\left(n^{-1/2}(2m-1)!!\,\frac{ f_3(\la)^n}{\la^{2m}}\right).
\end{equation}
\si

Let $d_1,\dots,d_n\ge 3$, meeting $\sum_i d_i=2m$, be such that there exists a graph
with the degree sequence $\bold d=(d_1,\dots,d_n)$; we call such $\bold d$ graphical.
Existence of  graphical $\bold d$'s for large $n,m$  is a weak consequence of  \eqref{CTheta}.
Let $g(\bold d)$ denote the total number of graphs for a graphical $\bold d$. Introduce
$G_{\bold d}$, a random graph distributed uniformly on the set of all $g(\bold d)$
graphs of a given graphical $\bold d$; obviously, $G_{\bold d}$ equals, in distribution,
$G^{(3)}(n,m)$, conditioned on $\{\bold d(G^{(3)}(n,m))=\bold d\}$. To handle $G_{\bold d}$
we use a random pairing model, Bollob\'as~\cite{eqref5}, defined as follows.

Introduce a partition of $[2m]$ into $n$ disjoint subsets
$Q_1,\dots,Q_n$, $|Q_i|=d_i$, and a set $\Omega$ of all $(2m-1)!!$ pairings $\omega$ of
$2m$ points in $[2m]$. Each $\omega$ induces a unique {\it multigraph\/}:
a pair $(u,v)\in\omega$ with $u,v\in Q_i$ becomes a loop at vertex $i$; a pair $(u,v)\in
\omega$ with $u\in Q_i$, $v\in Q_j$ becomes an edge $(i,j)$. Let $\omega$ be random,
distributed uniformly on $\Omega$. Let $MG_{\bold d}=MG_{\bold d}(\omega)$
denote the random multigraph induced by $\omega$. And let $\Omega_g(\bold d)$ be the set of
all graphical $\omega$'s, those for which $MG_{\bold d}(\omega)$ is a simple
graph, i.e. has neither loops nor multiple edges. Then $MG_{\bold d}(\omega)$, conditioned
on $ \omega\in \Omega_g(\bold d)$, coincides, in distribution, with $G_{\bold d}$. This implies
that, for any graph property ${\cal G}$,
$$
\pr(G_{\bold d}\in {\cal G})=\frac{\pr\bigl(\{MG_{\bold d}\in {\cal G} \}\cap \Omega_g(\bold d)
\bigr)}{\pr(\Omega_g(\bold d))}\le \frac{\pr\bigl(MG_{\bold d}\in {\cal G}\bigr)}{\pr(\Omega_g
(\bold d))}.
$$
Crucially,
\begin{equation}\label{gd}
g(\bold d)=\frac{(2m-1)!!}{\prod\limits_i d_i!}\,\pr(\Omega_g(\bold d)),
\end{equation}
\cite{eqref5}. We conclude that
\begin{align*}
\pr(G^{(3)}(n,m)\in {\cal G})=&\,C(n,m)^{-1}\sum_{\bold d}\pr(G_{\bold d}\in {\cal G})
g(\bold d)\\
\le&\, \frac{(2m-1)!!}{C(n,m)}\sum_{\bold d}\pr(MG_{\bold d}\in{\cal G})\,\prod_{i\in [n]}(1/d_i!);
\end{align*}
the sums are over all admissible graphical $\bold d$. So, by \eqref{CTheta}, uniformly
for all graph properties ${\cal G}$,
\begin{equation}\label{PcalG}
\pr(G^{(3)}(n,m)\in {\cal G})\le_b n^{1/2}\frac{\la^{2m}}{f_3(\la)^n}
\sum_{\bold d}\pr(MG_{\bold d}\in{\cal G})\,\prod_{i\in [n]}(1/d_i!).
\end{equation}
(For brevity, we write $A\le_b B$ when $A=O(B)$ uniformly over parameters involved,
and $B$ is too long to compose nicely with the big $O$ notation.)
This bound is perfectly tailored for ${\cal G}$'s implicit in \eqref{Pev} and \eqref{density=1}.
 \bi
{\bf Part (i):\/} Denote the probability in \eqref{Pev} by $P_{n1}$.
Suppose that for some $\omega$ and $A=A(\omega)\subseteq [n]$, $|A|\leq \eps_0\ln n$,
the sub(multi)graph of $MG_{\bold d}(\omega)$ induced by $A$ has more edges than vertices.
Then there exists $k=k(\omega)\in
[2,\eps_0\ln n]$ and point sets $A_{i_1},\dots, A_{i_k}$
such that the pairing $\omega$ contains $(k+1)$ pairs of points from $A_{i_1}\cup\cdots\cup A_{i_k}$.

Combining \eqref{PcalG} and the union bound, we have then
\begin{multline}\label{Pn1<}
P_{n1}\le_b n^{1/2}\frac{\la^{2m}}{f_3(\la)^n}\sum_{4\le k\leq\eps_0\ln n}\binom{n}{k}
\frac{\bigl(2(k+1)-1\bigr)!!\,\bigl(2(m-k-1)-1)!!}{(2m-1)!!}\\
\times\sum_{\bold d}\binom{d_{1:k}}{2(k+1)}\prod_{i=1}^n\frac{1}{d_i!},
\end{multline}
$d_{1:k}:=d_1+\cdots+d_k$. The second line sum in \eqref{Pn1<} is
\begin{multline*}
\sum_{3k\le d\le 2m}\binom{d}{2(k+1)}\sum_{d_1+\cdots+d_k=d\atop
d_i\ge 3}\prod_{i=1}^k\frac{1}{d_i!}\sum_{d_{k+1}+\cdots+d_n=2m-d\atop d_i\ge 3}
\prod_{i=k+1}^n\frac{1}{d_i!}\\
\le_b \sum_{3k\le d\le 2m}\binom{d}{2(k+1)}\,\frac{f_3(x)^k}{x^d}\,\frac{f_3(\la)^{n-k}}
{\la^{2m-d}},
\end{multline*}
for every $x>0$. (We have used a bound
\begin{equation}\label{d1d2..sum}
\sum_{\delta_1+\cdots+\delta_j=\delta\atop\delta_i\ge 3}\prod_{i=1}^j\frac{1}{\delta_i !}=
[y^{\delta}]\,f_3(y)^j\le\frac{f_3(y)^j}{y^{\delta}},\quad\forall\,y>0.)
\end{equation}
The ratio of two consecutive terms in the last sum is
$$
\frac{d+1}{d-2k-1}\,\frac{\la}{x}\leq\frac{3k+1}{k-1}\,\frac{\la}{x}\le 7\frac{\la}{x}=\frac{7}{8}<1,
$$
if
$$x=8\la.$$
So the sum is of order
$$
\binom{3k}{2(k+1)}\,\frac{f_3(8\la)^k}{(8\la)^{3k}}\,\frac{f_3(\la)^{n-k}}{\la^{2m-3k}}.
$$
Using this bound, $\binom{n}{k}\leq n^k/k!$, $2m/n=\la f_2(\la)/f_3(\la)$, and
$$
\bigl(2(k+1)-1\bigr)!!=\frac{\bigl(2(k+1)\bigr)!}{2^{k+1}(k+1)!},
$$
 we easily transform \eqref{Pn1<} into
 \begin{multline}\label{Pn1expl}
P_{n1}\le_b m^{-1/2}\sum_{k\le\eps_0\ln n}
\frac{k\,(3k)!}{2^k(k!)^3}\left[\frac{n}{2m}\,\frac{\la^3}{(8\la)^3}\,\frac{f_3(8\la)}{f_3(\la)}\right]^k\\
\le_b n^{-1/2}\sum_{k\le\eps_0\ln n}k\left[\frac{27}{2^{10}}\,\frac{f_3(8\la)}
{\la f_2(\la)}\right]^k\to 0,
\end{multline}
as
$$
\eps_0=(1/3)\ln^{-1}\left[\frac{27\, f_3(8\la)}{2^{10}\,\la f_2(\la)}\right].
$$
{\bf Part (ii):\/} Let $P_{n2}$ be the probability in \eqref{density=1}. This time we need to bound
the probability that there exists $A\subset [n]$, of cardinality $k\in
[\eps_0\ln n, n^{1-\rho_n}]$ that has
at least $\ell=\ell(k)=\lceil(1+\sigma_n)k\rceil$ edges. The counterpart of \eqref{Pn1<} is
\begin{multline}\label{Pn2<}
P_{n2}\le_b n^{1/2}\frac{\la^{2m}}{f_3(\la)^n}\sum_{\eps_0\ln n\le k\leq n^{1-\rho_n}}\binom{n}{k}
\frac{\bigl(2\ell-1\bigr)!!\,\bigl(2(m-\ell)-1)!!}{(2m-1)!!}\\
\times\sum_{\bold d}\binom{d_{1:k}}{2\ell}\prod_{i=1}^n\frac{1}{d_i!}.
\end{multline}
The bottom sum is bounded by a sum
$$
\sum_{3k\le d\le 2m}\binom{d}{2\ell}\,\frac{f_3(x)^k}{x^d}\,
\frac{f_3(\la)^{n-k}}{\la^{2m-d}},
$$
with the consecutive terms ratio bounded by
$3.01\la/x\le 0.76$, if $x=4\la$. So, like \eqref{Pn1expl},
\begin{equation}\label{Pn2expl}
\begin{aligned}
P_{n2}\le_b&\, n^{1/2}\sum_{\eps_0\ln n\le  k\leq  n^{1-\rho_n}}\binom{n}{k}
\frac{\bigl(2\ell-1\bigr)!!\,\bigl(2(m-\ell)-1)!!}{(2m-1)!!}\\
&\times \binom{3k}{2\ell}\gamma^k,\qquad
\gamma:=\,4^{-3}\frac{f_3(4\la)}{f_3(\la)}.
\end{aligned}
\end{equation}
Here
$$
\binom{3k}{2\ell}\le \binom{3k}{2k}\le_b \left(\frac{3^3}{2^2}\right)^k,
$$
and
$$
\frac{\bigl(2\ell-1\bigr)!!\,\bigl(2(m-\ell)-1)!!}{(2m-1)!!}=\frac{\binom{m}{\ell}}{\binom{2m}{2\ell}}\le\binom{m}{\ell}^{-1}.
$$
Using the last two bounds and
$$
\binom{n}{k}\le\left(\frac{en}{k}\right)^k,\quad m^{-1/2}\left(\frac{m}{\ell}\right)^{\ell}\le_b
\binom{m}{\ell},
$$
we simplify  \eqref{Pn2expl}  to
\begin{equation*}
P_{n2}\le_b\, n\sum_{\eps_0\ln n\le  k\leq
n^{1-\rho_n}}\left(\frac{n}{k}\right)^k\left(\frac{m}{\ell}\right)^{-\ell}\,\gamma_1^k,\quad
\gamma_1:=e\,3^3\,2^{-2}\,\gamma.
\end{equation*}
Here, since $\ell\ge (1+\sigma_n)k$,
\begin{align*}
\left(\frac{n}{k}\right)^k\left(\frac{m}{\ell}\right)^{-\ell}\gamma_1^k
\le&\,\left(\frac{n}{k}\right)^k\left(\frac{n}{k(1+\sigma_n)}\right)^{-k(1+\sigma_n)}\\
\le&\,
\left(\frac{n}{k}\right)^{-k\sigma_n}\cdot\bigl[(1+\sigma_n)^{-(1+\sigma_n)}\bigr]^k\gamma_1^k.
\end{align*}
The last expression is decreasing for $k\le n^{1-\rho_n}$, because its logarithmic derivative is
\begin{multline*}
-\sigma_n\ln\frac{n}{k}+\sigma_n-(1+\sigma_n)\ln(1+\sigma_n)+\ln\gamma_1\\
\le -\sigma_n\ln\frac{n}{k}+\ln\gamma_1\le -\sigma_n\rho_n\ln n +\ln\gamma_1\\
= -(\sigma_n\ln n)^{1/2}+\ln\gamma_1\to -\infty,
\end{multline*}
as $\sigma_n\ln n\to\infty$. So
$$
P_{n2}\le_b\exp\bigl[-\eps_0\sigma_n (\ln n)^2+O\bigl((\ln n)\ln\ln n\bigr)\bigr]\to 0,
$$
as $\sigma_n \ln n\gg \ln\ln n$.
\qed

\section{Moderately large, sparse Posa's  sets are unlikely.}

Let $d_{\text{max}}=d_{\text{max}}(n,m)$ denote the largest vertex degree in $G^{(3)}(n,m)$. Then
let $S$, $T$ be disjoint subsets of $[n]$, of cardinalities $s$ and $t$, with $t<2s$.
 In view of Lemma
 \ref{e<(1+eps)v}, part (i), we may and will confine ourselves to $s+t\ge\eps_0\ln n$.
Lemma \ref{Nu} asserts two upper bounds
 for the total number of subgraphs $G^{**}(S\cup T)$ with parameters $\bold s$,
$\bold t$ and $\boldsymbol{\xi}$. (See \eqref{xi's} for definition of $\xi_1$ and $\xi_2$.)
 Let us bound the number of ways to extend this subgraph to a  graph on $[n]$, of minimum degree $3$
at least, with $m$ edges.

Recall that the edge set of $G^{**}(S\cup T)$ does not contain edges between $T$-vertices.  So any such
extension of $G^{**}(S\cup T)$ is determined by an induced subgraph  $G(S^c)$. Let $d_i$ denote
the degree of vertex $i\in S^c$ in $G(S^c)$. An admissible $\bold d=\{d_i\}_{i\in S^c}$
meets the conditions
\begin{equation}\label{dgeq}
d_i\ge\left\{\alignedat2
& 3,\quad&&i\in S^c\setminus T,\\
&3-i,\quad&&i\in T_i,\,\,i=1,2,3,\endalignedat\right.
\end{equation}
and
\begin{equation}\label{sumd}
\sum_{i\in S^c}d_i =2m-2D,\quad D:=\mu+t_1=o(n).
\end{equation}
Then, by \eqref{gboldd} and the definition of $f_k(y)$,
the number of ways to extend a given  $G^{**}(S\cup T)$ is bounded above by
\begin{multline*}
\bigl(2(m-D)-1\bigr)!!\sum_{\bold d\text{ meets }\atop \eqref{dgeq}-\eqref{sumd}}
\prod_{i\in S^c}\frac{1}{d_i!}\\
=\,\bigl(2(m-D)-1\bigr)!!\,[y^{2(m-D)}]
\prod_{i=1}^3\!\!\left(\sum_{d\ge 3-i}\frac{y^d}{d!}\right)^{t_i}
\left(\sum_{d\ge 3}\frac{y^d}{d!}\right)^{n-s-t}\\
=\,\bigl(2(m-D)-1\bigr)!!\,[y^{2(m-D)}]\,\prod_{i=1}^3 f_{3-i}(y)^{t_i}\cdot f_3(y)^{n-s-t}.
\end{multline*}
By the Cauchy integral formula,
\begin{multline*}
[y^{2(m-D)}]\,\prod_{i=1}^3 f_{3-i}(y)^{t_i}\cdot f_3(y)^{n-s-t}\\
=\frac{1}{2\pi i}\oint\limits_{|y|=r}
\frac{1}{y^{2(m-D)+1}} \prod_{i=1}^3 f_{3-i}(y)^{t_i}\cdot f_3(y)^{n-s-t}\, dy.
\end{multline*}
Here $n-s-t\sim n$.  Using  $|f_k(y)|\le f_k(|y|)$, an inequality (\cite{refBP})
$$
|f_3(y)|\le f_3(|y|)\exp\left(-\frac{|y|-\text{Re }y}{4}\right),
$$
and selecting $r=\lambda$, we obtain
\begin{align*}
&\card{\frac{1}{2\pi i}\oint\limits_{|y|=r}
\frac{1}{y^{2(m-D)+1}} \prod_{i=1}^3 f_{3-i}(y)^{t_i}\cdot f_3(y)^{n-s-t}\, dy}\\
&\leq_b \frac{1}{\lambda^{2(m-D)}} \prod_{i=1}^3 f_{3-i}(\lambda)^{t_i}\cdot f_3(\lambda)^{n-s-t}
\int_{\th=-\pi}^\pi e^{-(n-o(n))\la(1-\cos\th)/4}d\th\\
&\leq_b \frac{1}{n^{1/2}\lambda^{2(m-D)}}
\prod_{i=1}^3 f_{3-i}(\lambda)^{t_i}\cdot f_3(\lambda)^{n-s-t}
\end{align*}

And so the number of extensions of a given  $G^{**}(S\cup T)$ is of order
 $$
N_{\text{ext}}(\bold s,\bold t,\boldsymbol\xi):=
\frac{\bigl(2(m-D)-1\bigr)!!}{n^{1/2}\lambda^{2(m-D)}} \prod_{i=1}^3
f_{3-i}(\lambda)^{t_i}\cdot f_3(\lambda)^{n-s-t},
$$
at most. Then,  multiplying $N_{\text{ext}}(\bold s,\bold t,\boldsymbol{\xi})$ by $N_1(\bold s,t,
\boldsymbol{\xi})$, the first bound given in
Lemma \ref{Nu}, we get an upper bound for the total
number of graphs on $[n]$ with $m$ edges, such that
$S\cup T$ induces a subgraph $G(S\cup T)$ with parameters $\bold s,\bold t,\boldsymbol
\mu$. Multiplying
$N(\bold s,t,\boldsymbol\xi)N_{\text{ext}}
(\bold s,\bold t,\boldsymbol\xi)$ by
$\binom{n}{s,t}\le n^{s+t}/s!t!$, and dividing by $C(n,m)$, the total
number of the $(n,m)$-graphs of minimum degree $3$ at least, we obtain a bound
$O\bigl(E_{n,m}(\bold s,\bold t,\boldsymbol\xi)\bigr)$
for the expected number of Posa's subgraphs with parameters $\bold s,\bold t,\boldsymbol\mu$, where
\begin{multline}\label{estmu}
E_{n,m}(\bold s,\bold t,\boldsymbol{\xi})= n^{s+t}\,\frac{(2\mu_1+\mu_2)!}{2^{\mu_1}\mu_1!}\,
\frac{\bigl(2(m-D)-1\bigr)!!}{(2m-1)!!}\exp\left[O\bigl(\sigma (s+t)\bigr)\right]\\
\times \frac{\lambda^{2D}}{2^{s_2+t_2}f_3(\lambda)^{s+t}}
\prod_{i=1}^3 f_{3-i}(\lambda)^{t_i}\cdot [f_3(1)]^{s_3+t_3}\\
\times\frac{1}{s!\,t!}\,t_1!\,\binom{s}{t_1}\binom{t}{t_1}\binom{s}{s_2}\binom{t-t_1}{t_2}.
\end{multline}
(See \eqref{mu=} and \eqref{mu1=} for $\mu_1,\mu_2$ expressed through $\xi_1$ and
$\xi_2$.) In view of \eqref{s3+t3}-\eqref{eS3T3} and
Lemma  \ref{e<(1+eps)v}, part (ii), if we allow only $s+t\le n^{1-\rho_n}$,
$\rho_n\to 0$, which we do, we need to
consider only $\bold s,\bold t,\boldsymbol{\xi}$ such that
\begin{equation}\label{stmu}
0\le s-t_1\le 2\sigma_n(s+t),\quad s_3+t_3\le 2\sigma_n(s+t),\quad \xi_1+\xi_2\le 2\sigma_n(s+t),
\end{equation}
where $\sigma_n\to 0$. (See Lemma \ref{e<(1+eps)v} for a more precise
definition of $\sigma_n,\rho_n$.) Our remaining task is to show that
the sum of $E_{n,m}(\bold s,\bold t,\boldsymbol{\xi})$ over
the admissible $(\bold s,\bold t,\boldsymbol{\xi})$ approaches zero.

To this end, let us first bound  $E_{n,m}(\bold s,\bold t,\boldsymbol{\xi})$ by a simpler
$E_{n,m}^*(\bold s,\bold t,\boldsymbol{\xi})$
times  $\exp\bigl(o(s+t)\bigr)$. First, by \eqref{sumd} and \eqref{stmu},
in the second line of \eqref{stmu}
\begin{align}
\frac{\lambda^{2D}}{2^{s_2+t_2}}\,f_0(\lambda)^{t_3}\,[f_3(1)]^{s_3+t_3}
=&\,\frac{\lambda^{2(s+t)}}{2^{t}}\exp\bigl(O(s_3+t_3+s-t_1)\bigr)\notag\\
=&\,\frac{\lambda^{2(s+t)}}{2^{t}}\exp\bigl(O(\sigma_n(s+t)\bigr)\label{la2D}.
\end{align}
Next, using
$$
(2a-1)!!=\frac{(2a)!}{2^a a!}=\Theta\left[\left(\frac{2a}{e}\right)^a\right],
$$
we obtain that the second fraction fraction in the first line of \eqref{estmu} is of order
\begin{align}
\left(\frac{e}{2m}\right)^D(1-D/m)^{m-D}=&\,(2m)^{-D}\,e^{O(D^2/m)}\notag\\
=&\,(2m)^{-t_1-\mu}\exp\bigl(O(n^{-\rho_n}(s+t))\bigr)\notag\\
=&(2m)^{-(s+t+\xi/2)}\exp\bigl(O(n^{-\rho_n}(s+t))\bigr)\label{2(m-D)},
\end{align}
$\xi:=\xi_1+\xi_2$. Further, by \eqref{st1}, \eqref{mu1,mu2} and \eqref{xi<},
\begin{equation}\label{2mu1+mu2}
\frac{(2\mu_1+\mu_2)!}{2^{\mu_1}\mu_1!}=\frac{(2s+\xi_1)!}{2^{2s-t}\bigl[s-t+t_1+(\xi_1-\xi_2)/2
\bigr]!}\exp\bigl(O(\sigma_n(s+t))\bigr).
\end{equation}
Given $\xi$, the last fraction attains its {\it maximum\/} at $\xi_1=\xi$,
$\xi_2=0$, and it is
\begin{multline}\label{1trin}
\frac{(2s+\xi)!}{2^{2s-t}(s-t+t_1+\xi/2)!}\\
=\frac{(s+t-t_1)!\,(\xi/2)!}{2^{2s-t}}
\binom{2s+\xi}{s-t+t_1+\xi/2,\,\,\xi/2,\,\,s+t-t_1}.
\end{multline}
The reason behind \eqref{1trin} is that the multinomial coefficients
are amenable to easy but sharp estimates. The factorial
$(s+t-t_1)!$ combined with $(t_1!/s!\,t!)\binom{s}{t_1}\binom{t}{t_1}$ in \eqref{estmu}
will later produce another  friendly trinomial coefficient.

Using an inequality
\begin{equation}\label{abc}
\binom{a+b+c}{a,b,c}\le\frac{(a+b+c)^{a+b+c}}{a^a\,b^b\,c^c},
\end{equation}
the trinomial coefficient in \eqref{1trin} is bounded above by $e^{H_1(s,\bold t,\xi)}$, where
\begin{multline}\label{H1}
H_1(\bold s,\bold t,\boldsymbol\xi)=(s-t+t_1+\xi/2)\ln\frac{2s+\xi}{s-t+t_1+\xi/2}\\
+\xi/2\ln\frac{2s+\xi}{\xi/2}+(s+t-t_1)\ln\frac{2s+\xi}{s+t-t_1}.
\end{multline}
Consider the first summand. Notice that
$$
s-t+t_1+\xi/2= 2s-t -(s-t_1)+\xi/2\le 2s-t +\xi/2,
$$
and $2s-t>0$. Suppose that $2s-t\ge \sigma_n^{1/2}(s+t)$. Then, as $s-t_1$ and $\xi$
are of order $O\bigl(\sigma_n(s+t)\bigr)$,  the summand is
$$
(2s-t)\ln\frac{2s}{2s-t}+O\bigl(\sigma_n^{1/2}(s+t)\bigr).
$$
If $2s-t\le  \sigma_n^{1/2}(s+t)$, then, as $x\ln (a/x)$ is increasing for $x\le a/e$, the
summand is bounded above crudely by
\begin{multline*}
(2s-t+\xi/2)\ln\frac{2s+\xi}{2s-t+\xi/2}\le 2\sigma_n^{1/2}(s+t)\ln\frac{3s}{\sigma_n^{1/2}(s+t)}\\
\le(2s-t)\ln\frac{2s}{2s-t}+2\sigma_n^{1/2}(\ln(1/\sigma_n) )(s+t).
\end{multline*}

Thus the summand is always
$$
(2s-t)\ln\frac{2s}{2s-t}+2\sigma_n^{1/2}(\ln(1/\sigma_n)) (s+t),
$$
at most. For the second summand in \eqref{H1},
$$
\xi/2\ln\frac{2s+\xi}{\xi/2}\le \sigma_n(s+t)\ln\frac{3s}{\sigma_n(s+t)}=
O\bigl(\sigma_n(\ln 1/\sigma_n)
(s+t)\bigr).
$$
The third summand in \eqref{H1} is
$$
t\ln\frac{2s}{t}+O\bigl(\sigma_n(s+t)).
$$
Thus
\begin{equation}\label{H1<}
H_1(\bold s,\bold t,\boldsymbol\xi)\le (2s-t)\ln \frac{2s}{2s-t}+t\ln\frac{2s}{t}+
O\bigl(\sigma_n^{1/2}(\ln(1/\sigma_n) )(s+t)
\bigr),
\end{equation}
uniformly for $2s-t>0$. The equation \eqref{2mu1+mu2} becomes
\begin{multline}\label{2mu1mu2better}
\frac{(2\mu_1+\mu_2)!}{2^{\mu_1}\mu_1!}\le \frac{(s+t-t_1)!\,(\xi/2)!}{2^{2s-t}}\\
\times\exp\left[(2s-t)\ln\frac{2s}{2s-t}+t\ln\frac{2s}{t}\right]\,
\cdot\exp\left[O\bigl(\sigma_n^{1/2}(\ln(1/\sigma_n) )(s+t)\right].
\end{multline}
Collecting  \eqref{la2D}, \eqref{2(m-D)} and \eqref{2mu1mu2better}, we conclude that
$$
E_{n,m}(\bold s,\bold t,\boldsymbol{\xi})\le E_{n,m}^*(\bold s,\bold t,\boldsymbol\xi)
\exp\left[O\bigl(\sigma_n^{1/2}(\ln(1/\sigma_n) )(s+t)\right],
$$
where
\begin{align}
E_{n,m}^*(\bold s,\bold t,\boldsymbol\xi)=&\,\frac{(\xi/2)!}{m^{\xi/2}}\,
\frac{n^{s+t}\lambda^{2(s+t)}
f_2(\la)^sf_1(\la)^{t-t_1}}
{(2m)^{s+t}f_3(\la)^{s+t}\,2^{2s}}\notag\\
&\times \exp\left[(2s-t)\ln\frac{2s}{2s-t}+t\ln\frac{2s}{t}\right]\notag\\
&\times\frac{(s+t-t_1)!}{s!\,t!}\,
\,t_1!\,\binom{s}{t_1}\binom{t}{t_1}\binom{s}{s_2}\binom{t-t_1}{t_2}.\label{E*(bs,bt,xi)}
\end{align}
Here,  recalling again \eqref{stmu},
$$
s_3,\,t_3\le 2\sigma_n(s+t),\quad \xi\le 2\sigma_n(s+t).
$$
Subject to this constraint,  let us bound
$\sum_{s_2,t_2,\xi}E^{*}_{n,m}(\bold s,\bold t,\boldsymbol\xi)$. First of all,
\begin{align*}
\sum_{\xi}\frac{(\xi/2)!}{m^{\xi/2}}\le&\, \sum_{\xi}\left(\frac{e\xi}{2m}\right)^{\xi/2}\\
\le&\,\sum_{\xi\ge 0}\left(\frac{e\sigma_n(s+t)}{m}\right)^{\xi/2}\to 1.
\end{align*}
Secondly,
\begin{align*}
\sum_{s_2+s_3=s\atop s_3\le 2\sigma_n(s+t)}\binom{s}{s_2}=&\,\sum_{s_3\le
2\sigma_n(s+t)}\binom{s}{s_3}\\
\leq_b&\,\binom{s}{2\sigma_n(s+t)}\le\left(\frac{es}{2\sigma_n(s+t)}\right)^{2\sigma_n(s+t)}\\
=&\,\exp\left[O\bigl(\sigma_n(\ln(2/\sigma_n))(s+t)\bigr)\right].
\end{align*}
Likewise
$$
\sum_{t_2+t_3=t-t_1\atop t_3\le 2\sigma_n(s+t)}\binom{t-t_1}{t_2}=
\exp\left[O\bigl(\sigma_n(\ln(2/\sigma_n))(s+t)\bigr)\right].
$$
Observing also that
$$
\frac{(s+t-t_1)!}{s!\,t!}\,\,t_1!\,\binom{s}{t_1}\binom{t}{t_1}=\binom{s+t-t_1}{s-t_1,\,t_1,\,t-t_1},
$$
we then have
$$
\sum_{s_2,t_2,\xi}E^{*}_{n,m}(\bold s,\bold t,\xi)\le  E_{n,m}(s,t,t_1)
\exp\left[O\bigl(\sigma_n(\ln(2/\sigma_n))(s+t)\bigr)\right],
$$
where
\begin{multline}\label{Enmstt1}
E_{n,m}(s,t,t_1):=\,\frac{n^{s+t}\lambda^{2(s+t)}f_2(\la)^sf_1(\la)^{t-t_1}}
{(2m)^{s+t}f_3(\la)^{s+t}\,2^{2s}}\\
\times\exp\left[(2s-t)\ln\frac{2s}{2s-t}+t\ln\frac{2s}{t}\right]
\cdot \binom{s+t-t_1}{s-t_1,t_1,t-t_1}.
\end{multline}
The trinomial coefficient in \eqref{Enmstt1} is bounded above by
\begin{equation}\label{trin}
\exp\left[(s-t_1)\ln\frac{s+t-t_1}{s-t_1}+t_1\ln\frac{s+t-t_1}{t_1}+(t-t_1)\ln\frac{s+t-t_1}{t-t_1}\right].
\end{equation}
Recall that
\begin{equation}\label{recall}
s-2\sigma_n(s+t)\le t_1 \le s\Longrightarrow 0\le s-t_1\le 2\sigma_n(s+t).
\end{equation}
Since
$$
(s-t_1)\ln\frac{s+t-t_1}{s-t_1}= (s-t_1)\ln\frac{t}{s-t_1}+O\bigl((s-t_1)^2/t\bigr),
$$
and $x\ln (t/x)$ is increasing for $x<t/e$, we obtain
\begin{equation}\label{s-t1}
\begin{aligned}
(s-t_1)\ln\frac{s+t-t_1}{s-t_1}\leq&\,
2\sigma_n(s+t)\ln\frac{t}{2\sigma_n(s+t)}+O\bigl(\sigma_n^2(s+t)
\bigr)\\
=&\,O\bigl(\sigma_n\ln(1/\sigma_n)(s+t)\bigr),
\end{aligned}
\end{equation}
where $\sigma_n\ln(1/\sigma_n)\to 0$, as $\sigma_n\to 0$. Furthermore,
\begin{equation}\label{t1}
t_1\ln\frac{s+t-t_1}{t_1}=s\ln\frac{t}{s}+O(\sigma_n(s+t)).
\end{equation}
Turn to the last summand in \eqref{trin}.  By \eqref{recall},
$$
0\le t-t_1\le t-s +2\sigma_n(s+t)<  t-s +3\sigma_n(s+t).
$$
Further
\begin{align*}
(t-t_1)\ln\frac{s+t-t_1}{t-t_1}=&\,(t-t_1)\ln\frac{t}{t-t_1}+O(s-t_1)\\
=&\,(t-t_1)\ln\frac{t}{t-t_1}+O(\sigma_n(s+t)).
\end{align*}
Suppose that
\begin{equation}\label{sun}
t-s +3\sigma_n(s+t)<t/e
\end{equation}
which is equivalent to
$$t<\frac{s(1-3\sigma_n)}{1-e^{-1}+3\sigma_n}.$$
Then, for $t>s$,
\begin{equation}\label{t>s,t-t1}
\begin{aligned}
(t-t_1)\ln\frac{t}{t-t_1}\le& \bigl(t-s+3\sigma_n(s+t)\bigr)\ln\frac{t}{t-s+3\sigma_n(s+t)}\\
\le&(t-s)\ln\frac{t}{t-s}+O\left[\sigma_n(s+t)\ln\frac{t}{\sigma_n(s+t)}\right]\\
=&\,(t-s)\ln\frac{t}{t-s}+O(\sigma_n\ln(1/\sigma_n) (s+t)).
\end{aligned}
\end{equation}
If  $t\le s$, then \eqref{sun}
\begin{equation}\label{t<s,t-t1}
\begin{aligned}
(t-t_1)\ln\frac{t}{t-t_1}\le&\, \bigl(t-s+3\sigma_n(s+t)\bigr)\ln\frac{t}{t-s+3\sigma_n(s+t)}\\
\le&\,3\sigma_n(s+t)\ln\frac{t}{3\sigma_n(s+t)}\\
=&\,O(\sigma_n\ln(1/\sigma_n) (s+t)).
\end{aligned}
\end{equation}
Suppose that
$$
t\ge\frac{s(1-3\sigma_n)}{1-e^{-1}+3\sigma_n}.
$$
Then $t-s=\Theta(s)$, and so
\begin{equation}\label{t>}
(t-t_1)\ln\frac{t}{t-t_1}-(t-s)\ln\frac{t}{t-s}=O\bigl(s-t_1\bigr)=O\bigl(\sigma_n(s+t)\bigr).
\end{equation}
Combining \eqref{Enmstt1}-\eqref{t>}, we obtain
$$
E_{n,m}(s,t):=\,\sum_{t_1}E_{n,m}(s,t,t_1)\leq E_{n,m}^*(s,t)
\exp\left[O\bigl(\sigma_n(\ln 1/\sigma_n)(s+t)\bigr)\right],
$$
where
\begin{multline}\label{Enm*st}
E_{n,m}^*(s,t):=\,t\,\frac{n^{s+t}\lambda^{2(s+t)}}{(2m)^{s+t}f_3(\la)^{s+t}\,2^{2s}}\,f_2(\la)^sf_1(\la)^{t-s}\\
\times\,\exp\left[(2s-t) \ln\frac{2s}{2s-t}+t\ln\frac{2s}{t}+
s\ln\frac{t}{s}+(t-s)^+\ln\frac{t}{(t-s)^+}\right],
\end{multline}
where we define $x^+=\max\{0,x\}$, and $0\ln (t/0)=0$.
\si

The rest is a bit of calculus. Recalling that $2m/n=\la f_2(\la)/f_3(\la)$,
and setting $t=xs$, we write
$$
E_{n,m}^*(s,t)=sxe^{sH_1(x)},
$$
where
\begin{multline*}
H_1(x)=\,(1+x)\ln\la-x\ln f_2(\la)+(x-1)\ln f_1(\la)\\
+(2-x)\ln\frac{1}{2-x}+(x-1)\ln\frac1x+(x-1)^+\ln\frac{ x}{(x-1)^+}.
\end{multline*}
The exponent in \eqref{Enm*st} is
\begin{multline*}
(2s-t)\ln \frac{2s}{2s-t}+ t\ln\frac{2s}{t}+s\ln \frac{t}{s}\\
=2s\ln 2 +(2s-t)\ln \frac{s}{2s-t}+ t\ln\frac{s}{t}+s\ln \frac{t}{s}\\
=2s\ln 2 +(2s-t)\ln \frac{s}{2s-t}+(t-s)\ln\frac{s}{t},
\end{multline*}
and the term $2s\ln 2$ cancels $2^{2s}$ in the first line fraction denominator.
Since
$$
H^\prime_1(x)=\left\{\begin{aligned}
&\ln\left(\frac{\la f_1(\la)}{f_2(\la)}\frac{2-x}{x}\right)+\frac{1}{x},\qquad\, \,x<1,\\
&\ln\left(\frac{\la f_1(\la)}{f_2(\la)}\,\frac{2-x}{x-1}\right),\quad x\in (1,2),\end{aligned}\right.
$$
and $\la f_1(\la)/f_2(\la)>2$, we see that $H_1(x)$ is unimodal on $(0,2)$,
and attains its maximum at $x^*\in (1,2)$
$$
 x^*=\frac{1+2\la f_1(\la)/f_2(\la)}{1+\la f_1(\la)/f_2(\la)},
$$
and
\begin{equation}\label{Hx*}
H_1(x^*)=\ln\left[\frac{\la^2}{f_2(\la)}\,\bigl(1+\la f_1(\la)/f_2(\la)\bigr)\right]
\end{equation}
 Maple shows that the function on the RHS of \eqref{Hx*} increases with $\la$ and it is zero at
$\lambda^*=5.162717...$. At the first glance it would seem necessary to put a constraint
$\la>\la^*$in order to claim that , for those $\la$'s, whp there are no Posa's sets of
cardinality $|S|+|T|\le n^{1-o(1)}$.
\si

We can do better though! Indeed, by the unimodality of $H_1(x)$,
\begin{multline}\label{maxH1}
\max\bigl\{H_1(x)\,:\,x\in [0,1+\sigma_n^{-1/2}]\cup [2-\sigma_n^{1/2},2]\bigr\}\\
=\max\bigl\{H_1(1+\sigma_n^{-1/2}),\,H_1(2-\sigma_n^{1/2})\bigr\}\\
=\max\left\{\ln\frac{\la^2}{f_2(\la)},\,\ln\frac{\la^3 f_1(\la)}{f_2(\la)^2}\right\}+
O\bigl(\sigma_n^{1/2}\ln(1/\sigma_n)\bigr)\\
=\ln\frac{\la^3 f_1(\la)}{f_2(\la)^2}+O\bigl(\sigma_n^{1/2}\ln(1/\sigma_n)\bigr),
\end{multline}
as $\la f_1(\la)/f_2(\la)>2$. As for $x=t/s\in [1+\sigma_n^{1/2},2-\sigma_n^{1/2}]$, we use
\eqref{2Nu<} instead of \eqref{1Nu<} and improve  the bound \eqref{Enm*st} by the factor
$$
(s+t)^2\exp\left[-(2s-t)\ln\frac{s}{2s-t}-(t-s)\ln\frac{s}{t-s}.
\right]
$$
So we can re-define
$$
E_{n,m}^*(s,t)=s^3x(1+x)^2e^{sH_2(x)},
$$
where
$$
H_2(x)=(1+x)\ln\la -x\ln f_2(\la)+(x-1)\ln f_1(\la),
$$
a linear function! Now
$$
\max\{H_2(x)\,:\,x\in [1,2]\} = H_2(2)=\ln\frac{\la^3 f_1(\la)}{f_2(\la)^2},
$$
and this function decreases with $\la$.

Indeed, introducing $F(\la)=\la/f_1(\la)$ that decreases from $1$ at
$0+$ to $0$ at $\infty$, we have
$$
\frac{\la^3 f_1(\la)}{f_2(\la)^2}=\frac{\la^3f_1(\la)}{\bigl(f_1(\la)-\la\bigr)^2}=
\la^2\frac{F(\la)}{\bigl(1-F(\la)\bigr)^2}.
$$
So
\begin{align*}
\frac{d}{d\la}\frac{\la^3 f_1(\la)}{f_2(\la)^2}=&\,
2\la\frac{F(\la)}{\bigl(1-F(\la)\bigr)^2}+\la^2\frac{1+F(\la)}{\bigl(1-F(\la)\bigr)^3}\,
F^\prime(\la)\\
&\left(\text{using }F^\prime(\la)=\frac{1}{f_1(\la)}-
\frac{\la e^{\la}}{f_1(\la)^2}=\la^{-1}\bigl(F(\la)-
e^{\la}F(\la)^2\bigr)\right)\\
=&\,\frac{\la F(\la)}{\bigl(1-F(\la)\bigr)^3}
\left[3-F(\la)-e^{\la}F(\la)\bigl(1+F(\la)\bigr)\right]\\
=&\,\frac{\la F(\la)}{\bigl(1-F(\la)\bigr)^3(e^{\la}-1)}\,D(\la);
\end{align*}
here
$$
D(\la)=(3-\la)e^{2\la}-(6+\la^2)e^{\la}+\la+3=\sum_{j\ge 4}d_j\la^j,
$$
and
$$
d_j=3\cdot 2^j-j2^{j-1}-j(j-1)-6,\quad j\ge 4.
$$
By induction on $j\ge 4$, $d_j<0$ for all $j\ge 4$. Hence
$$
\frac{d}{d\la}\frac{\la^3 f_1(\la)}{f_2(\la)^2}<0,\quad\forall\,\la>0.
$$

Maple shows that $\la^3f_1(\la)/f_2(\la)^2$ attains value $1$ at
$$
\la^{**}=4.789771...
$$
The corresponding average vertex degree
$$
c^{**}=\frac{\la^{**}f_2(\la^{**})}{f_3(\la^{**})}=5.323132...
$$
\bi
{\bf It follows that for $m\ge 2.662 n$ the expected number of the
likely Posa's sets $(S,T)$ of size $|S|+|T|\le
n^{1-o(1)}$ approaches zero as $n,m\to\infty$.\/}

\begin{Remark}
As a final remark, observe that
within the constraints on the $G(S\cup T)$,
the dominant contribution to the total number of sparse Posa sets
$(S,T)$ comes from $G(S\cup T)$ very close to an alternating cycle
on $S + T_2, |T_2|=|S|=s$, with the $s$ pendant $T_1$ vertices attached to
to $S$-vertices.  It is not difficult to get directly the asymptotic expected number of
such extreme subgraphs  in our random graph, and it turns out essentially the same
as the current estimate.

What this likely means is that it is fruitless to search for another constraint
on $G(S\cup T)$ with a potential to further decrease  $\lambda^{**}$ via a sharper
bound for the expected number of Posa sets $(S,T)$.

\end{Remark}


\begin{thebibliography}{99}

\bibitem{AKS} M. Ajtai, J. Koml\'os and E. Szemer\'edi,
First occurrence of Hamilton cycles in random graphs,
\textit{Cycles in graphs (Burnaby, B.C.)} (1982) 173-178.

\bibitem{FP1} J. Aronson, A.M. Frieze and B. Pittel, Maximum matchings in sparse random
graphs: Karp-Sipser revisited,
\textit{Random Structures and Algorithms} \textbf{12} (1998) 111-178.







\bibitem{BF} T. Bohman and A.M. Frieze, Hamilton cycles in 3-out,
\textit{Random Structures and Algorithms},  \textbf{35} (2009) 393-417.

\bibitem{eqref5} B. Bollob\'as, A probabilistic proof of an asymptotic formula for the
number of labelled regular graphs, \textit{European J. Comb.}, \textbf{1}, (1980) 311--316.

\bibitem{eqref6} B. Bollob\'as, The evolution of random graphs,
\textit{Trans. Amer. Math. Soc.}, \textbf{286}, (1984) 257--274.

\bibitem{eqref7} B. Bollob\'as, {Random Graphs, Second Edition}
\textit{Cambridge University Press},  (2001).

\bibitem{Bo2} B. Bollob\'as, Almost all regular graphs are Hamiltonian,
\textit{European Journal on Combinatorics} \textbf{4} (1983) 97-106.

\bibitem{BCFF} B. Bollob\'as, C. Cooper, T. Fenner and A.M. Frieze,
On Hamilton cycles in sparse random graphs with minimum degree at least $k$,
\textit{Journal of Graph Theory} \textbf{34} (2000) 42-59.

\bibitem{BFF} B. Bollob\'as, T. Fenner and A.M. Frieze,
Hamilton cycles in random graphs with minimal degree at least $k$,
\textit{in A tribute to Paul Erdos, edited by A.Baker, B.Bollobas and A.Hajnal} (1990) 59-96.


\bibitem{eqref9} P. Erd\H os and A. R\'enyi,  On the evolution of random graphs,
\textit{Publ. Math. Hungar. Acad. Sci.}, \textbf{5},  (1960) 17--61.

\bibitem{FF} T. Fenner and A.M. Frieze, Hamiltonian cycles in random regular graphs,
\textit{Journal of Combinatorial Theory B} \textbf{40} (1984) 103-112.

\bibitem{dlV} W. Fernandez de la Vega, Long paths in sparse random graphs,
\textit{Studia Sci. Math. Hungar.} \textbf{14} (14) (1979) 335-340.

\bibitem{F} A.M. Frieze, Finding hamilton cycles in sparse random graphs,
\textit{Journal of Combinatorial Theory B}, \textbf{44} (1988) 230-250.

\bibitem{F2G} A.M. Frieze, On Hamilton Cycles in Random Graphs with Minimum Degree at Least Three.

\bibitem{FP2} A.M. Frieze and B. Pittel,
Perfect matchings in random graphs with prescribed minimal degree,
\textit{Trends in Mathematics, Birkhauser Verlag, Basel}, (2004) 95-132.




\bibitem{KaSi} R.M. Karp and M. Sipser, Maximum matchings in sparse random
graphs, \textit{Proceedings of the 22nd Annual IEEE Symposium on Foundations of
Computing} (1981) 364-375.

\bibitem{KS1} J. Koml\'os and E. Szemer\'edi, Hamilton cycles in random graphs,
\textit{Infinite and finite sets (Colloq., Keszthely, 1973;
dedicated to P. Erdos on his 60th birthday)},
(1975) 1003-1010.

\bibitem{Kor} A. Korshunov,  A solution of a
problem of P. Erd\H{o}s and A. R\'enyi about Hamilton cycles in non-oriented
graphs, \textit{Metody Diskr. Anal. Teoriy Upr. Syst. Sh. Trudov.}
\textbf{31} (1977) 17-56 (in Russian).

\bibitem{KS2} J. Koml\'os and E. Szemer\'edi,
Limit distribution for the existence of Hamiltonian cycles in a random graph,
\textit{Discrete Mathematics} \textbf{43} (1983) 55-63.









\bibitem{refBP} B. Pittel, On tree census and the giant component in sparse random graphs,
\textit{Random Structures and Algorithms} \textbf{1} (1990) 311-342.


\bibitem{eqref21} B. Pittel and N. C. Wormald, Counting connected graphs inside-out,
\textit{J. Comb. Theory, Ser. B}, \textbf{93}, (2005) 122--172.



\bibitem{Po} L. P\'osa, Hamiltonian circuits in random graphs,
\textit{Discrete Mathematics} \textbf{14} 359-364.

\bibitem{eqreRW} R. Robinson and N. Wormald, Almost all regular graphs are Hamiltonian,
\textit{Random Structures and Algorithms} \textbf{5}m (1994) 363-374.





\end{thebibliography}
\end{document}